\documentclass[]{article}

\usepackage{amssymb,amsmath,hyperref}
\usepackage{color,graphics,graphicx}
\usepackage{color,cite}
\usepackage{amsthm}	% Fancy theorem-type environments
\usepackage[margin=1in]{geometry}
\newtheorem{thm}{Theorem}
\newtheorem{asmpt}{Assumption}
\newtheorem{rmrk}{Remark}
\title{Adaptive Differential Thrust Methodology for Lateral/Directional Stability of an Aircraft with a Completely Damaged Vertical Stabilizer}

\author{
  Long Lu%
    \thanks{Graduate Student, Department of Aerospace Engineering, AIAA Student Member, long.lu@sjsu.edu.}
\  and \ Kamran Turkoglu\thanks{Corresponding Author, Assistant Professor, Department of Aerospace Engineering, AIAA Member, kamran.turkoglu@sjsu.edu}\\
  {\normalsize\itshape San Jos\'{e} State University, San Jose, CA 95192, USA} \\ \\
 }
\begin{document}
\maketitle
%\noindent \emph{Keywords:} MRAS Adaptive control, Lyapunov stability, Directional stability, Flight control, Flight dynamics
\begin{abstract}

This paper investigates the utilization of differential thrust to help a commercial aircraft with a damaged vertical stabilizer regain its lateral/directional stability. In the event of an aircraft losing its vertical stabilizer, the consequential loss of the lateral/directional stability and control is likely to cause a fatal crash. In this paper, an aircraft with a completely damaged vertical stabilizer is investigated, and a unique differential thrust based adaptive control approach is proposed to achieve a stable flight envelope. The propulsion dynamics of the aircraft is modeled as a system of differential equations with engine time constant and time delay terms to study the engine response time with respect to a differential thrust input. The proposed differential thrust control module is then presented to map the rudder input to differential thrust input. Model reference adaptive control based on the Lyapunov stability approach is implemented to test the ability of the damaged aircraft to track the model aircraft's (reference) response in an extreme scenario. Investigation results demonstrate successful application of such differential thrust approach to regain lateral/directional stability of a damaged aircraft with no vertical stabilizer. Finally, the conducted robustness and uncertainty analysis results conclude that the stability and performance of the damaged aircraft remain within desirable limits, and demonstrate a safe flight mission through the proposed adaptive control methodology.

\end{abstract}

\section*{Nomenclature}

\begin{tabbing}
  XXXXX \= \kill% this line sets tab stop
$A/C$ 		\> aircraft \\
$A$ 		\> state matrix\\
%$A_d$ 		\> state matrix of the damaged aircraft \\
%$A_m$ 		\> state matrix of the model aircraft \\
%$A_n$ 		\> state matrix of the nominal (undamaged) aircraft \\ 
$B$ 		\> input matrix\\
%$B_d$ 		\> input matrix of the damaged aircraft \\
%$B_m$ 		\> input matrix of the model aircraft \\
%$B_n$ 		\> input matrix of the nominal (undamaged) aircraft \\  
$b$ 		\> aircraft wing span \\
$C$ 		\> output matrix\\ 
%$C_d$ 		\> output matrix of the damaged aircraft\\
%$C_n$ 		\> output matrix of the nominal aircraft\\
${C_L}_{i}$		\> dimensionless derivative of rolling moment, $i=p, r, \beta, \delta_a, \delta_r$\\
%${C_L}_p$		\> dimensionless derivative of rolling moment coefficient with respect to roll rate \\
%${C_L}_r$		\> dimensionless derivative of rolling moment coefficient with respect to yaw rate \\
%${C_L}_{\beta}$		\> dimensionless derivative of rolling moment coefficient with respect to side slip angle \\
%${C_L}_{\delta a}$		\> dimensionless derivative of rolling moment coefficient with respect to aileron deflection \\
%${C_L}_{\delta r}$		\>  dimensionless derivative of rolling moment coefficient with respect to rudder deflection \\
%${{C_L}_\alpha}_v$		\> lift curve slope of the vertical stabilizer \\
${C_N}_{i}$		\> dimensionless derivative of yawing moment, $i=p, r, \beta, \delta_a, \delta_r$\\
%${C_N}_p$		\> dimensionless derivative of yawing moment coefficient with respect to roll rate \\
%${C_N}_r$		\> dimensionless derivative of yawing moment coefficient with respect to yaw rate \\
%${C_N}_{\beta}$		\> dimensionless derivative of yawing moment coefficient with respect to side slip angle \\
%${C_N}_{\delta a}$		\> dimensionless derivative of yawing moment coefficient with respect to aileron deflection \\ 
%${C_N}_{\delta r}$		\> dimensionless derivative of yawing moment coefficient with respect to rudder deflection \\ 
${C_Y}_{i}$		\> dimensionless derivative of side force, $i=p, r, \beta, \delta_a, \delta_r$\\
%${C_Y}_p$		\> dimensionless derivative of side force coefficient with respect to roll rate \\
%${C_Y}_r$		\> dimensionless derivative of side force coefficient with respect to yaw rate \\
%${C_Y}_{\beta}$		\> dimensionless derivative of side force coefficient with respect to side slip angle \\
%${C_Y}_{\delta a}$		\> dimensionless derivative of side force coefficient with respect to aileron deflection \\
%${C_Y}_{\delta r}$		\> dimensionless derivative of side force coefficient with respect to rudder deflection \\ 
$D$ 		\> state transition matrix\\
%$D_d$ 		\> state transition matrix of the damaged aircraft \\
%$D_n$ 		\> state transition matrix of the nominal aircraft\\
$\frac{d\sigma}{d\beta}$ 		\> change in side wash angle with respect to change in side slip angle\\
$g$ 		\> gravitational acceleration \\
%$h$ 		\> altitude \\
$\overline{I_{xx}}$ 		\> normalized mass moment of inertia about the x axis \\
$\overline{I_{xz}}$ 		\> normalized product of inertia about the xz axis \\
$\overline{I_{zz}}$ 		\> normalized mass moment of inertia about the z axis \\
$L_{i}$		\> dimensional derivative of rolling moment, $i=p, r, \beta, \delta_a, \delta_r$\\
%$L_p$		\> dimensional derivative of rolling moment coefficient with respect to roll rate \\
%$L_r$		\> dimensional derivative of rolling moment coefficient with respect to yaw rate \\
%$L_{\beta}$		\> dimensional derivative of rolling moment coefficient with respect to side slip angle \\
%$L_{\delta a}$		\> dimensional derivative of rolling moment coefficient with respect to aileron deflection \\
%$L_{\delta r}$		\>  dimensional derivative of rolling moment coefficient with respect to rudder deflection \\
%$M_u$		\> dimensional derivative of pitching moment coefficient with respect to linear velocity in the x-direction \\
%$M_q$		\> dimensional derivative of pitching moment coefficient with respect to pitch rate \\
%$M_{\dot{w}}$		\> dimensional derivative of pitching moment coefficient with respect to acceleration in the z-direction \\
$L_v$ 		\>  vertical stabilizer lift force\\
$l_v$ 		\> distance from the vertical stabilizer aerodynamic center to the aircraft center of gravity\\
$MRAS$ 		\> Model Reference Adaptive System \\
$m$ 		\> aircraft mass \\
$N_{i}$		\> dimensional derivative of yawing moment, $i=p, r, \beta, \delta_a, \delta_r, \delta T$\\
%$N_p$		\> dimensional derivative of yawing moment coefficient with respect to roll rate \\
%$N_r$		\> dimensional derivative of yawing moment coefficient with respect to yaw rate \\
%$N_{\beta}$		\> dimensional derivative of yawing moment coefficient with respect to side slip angle \\
%$N_{\delta a}$		\> dimensional derivative of yawing moment coefficient with respect to aileron deflection \\
%$N_{\delta r}$		\> dimensional derivative of yawing moment coefficient with respect to rudder deflection \\
%$N_{\delta T}$		\> dimensional derivative of yawing moment coefficient with respect to differential thrust \\
$p$ 		\> roll rate \\
$r$ 		\> yaw rate \\
$S$ 		\> aircraft wing area \\
$S_v$ 		\> vertical stabilizer area \\
$T$ 		\> engine thrust \\
$T_c$ 		\> engine thrust command \\
$t$		\> time \\
$t_d$		\> time delay\\
%$U$ 		\> velocity in the x direction \\
%$V$ 		\> velocity in the y direction \\
$V_v$ 		\> vertical stabilizer volume ratio\\
$v$ 		\> airspeed \\
%$W$ 		\> velocity in the z direction \\
$W$ 		\> aircraft weight \\
$Y_{i}$		\> dimensional derivative of side force, $i=p, r, \beta, \delta_a, \delta_r$\\
%$Y_p$		\> dimensional derivative of side force coefficient with respect to roll rate \\
%$Y_r$		\> dimensional derivative of side force coefficient with respect to yaw rate \\
%$Y_{\beta}$		\> dimensional derivative of side force coefficient with respect to side slip angle \\
%$Y_{\delta a}$		\> dimensional derivative of side force coefficient with respect to aileron deflection \\
%$Y_{\delta r}$		\> dimensional derivative of side force coefficient with respect to rudder deflection \\
$y_e$ 		\> distance from the outermost engine to the aircraft center of gravity\\
$z_v$ 		\> distance from the vertical stabilizer center of pressure to the fuselage center line\\
$\alpha$ 		\> angle of attack \\
$\beta$ 		\> side slip angle \\
$\gamma$  	\> flight path angle \\
$\delta a$		\> aileron deflection \\
$\delta r$		\> rudder deflection \\
$\Delta T$		\> collective thrust \\
$\delta T$		\> differential thrust \\
$\zeta$		\> damping ratio \\
$\eta$  	\> efficiency factor \\
$\theta$  	\> pitch angle \\
$\rho$		\> air density \\
$\tau$   	\> time constant \\
$\phi$    	\> roll angle \\
$\omega$		\> bandwidth frequency \\
$\dot{( \ )}$ 	\>  first order time derivative \\
$\ddot{( \ )}$ 	\>  second order time derivative \\
%$( \ )_E$	\>  with respect to Earth component \\
%$( \ )_e$	\>  equilibrium component \\
$\bar{( \ )}$ 	\>  trimmed value \\
$( \ )_d$	\>  damaged aircraft component \\
$( \ )_m$	\>  model aircraft component \\
$( \ )_n$	\>  nominal (undamaged) aircraft component \\
 \end{tabbing}

\section{Introduction} \label{Introduction} 

The vertical stabilizer is an essential element in providing the aircraft with its directional stability characteristic while ailerons and rudder serve as the primary control surfaces of the yawing and banking maneuvers. In the event of an aircraft losing its vertical stabilizer, the sustained damage will cause lateral/directional stability to be compromised, and the lack of control is likely to result in a fatal crash. Notable examples of such a scenario are the crash of the American Airline 587 in 2001 when an Airbus A300-600 lost its vertical stabilizer in wake turbulence, killing all passengers and crew members \cite{NTSB_flight587} and the crash of Japan Airlines Flight 123 in 1985 when a Boeing 747-SR100 lost its vertical stabilizer leading to an uncontrollable aircraft, resulting in 520 casualties \cite{faa_flight123}.

However, not all situations of losing the vertical stabilizer resulted in a total disaster. In one of those cases, the United Airlines Flight 232 in 1989 \cite{NTSB_flight_232}, differential thrust was proved to be able to make the aircraft controllable. Another remarkable endeavor is the landing of the Boeing 52-H even though the aircraft lost most of its vertical stabilizer in 1964 \cite{Hartnett06}.

Research on this topic has been conducted with two main goals: to understand the response characteristics of the damaged aircraft such as the work of Bacon and Gregory \cite{Bacon_Gregory_07}, Nguyen and Stepanyan \cite{NguyenStepanyan10}, and Shah \cite{Shah08}, as well as to come up with an  automatic control algorithm to save the aircraft from disasters, where similar such as the work of Burcham et al. \cite{Burcham_et_al_98}, Guo et al. \cite{Guo_et_al_11}, Liu et al. \cite{Liu_et_al_09}, Tao and Ioanou \cite{Tao_Ioannou_88}, and Urnes and Nielsen \cite{UrnesNielsen10}.

Notable research on the topic of a damaged transport aircraft includes the work of Shah \cite{Shah08}, where a wind tunnel study was performed to evaluate the aerodynamic effects of damages to lifting and stability/control surfaces of a commercial transport aircraft. In his work, Shah \cite{Shah08} studied this phenomenon in the form of partial or total loss of the wing, horizontal, or vertical stabilizers for the development of flight control systems to recover the damaged aircraft from adverse events. 

In literature, there exists a similar study conducted by Stepanyan et al. \cite{Stepanyan09}, which provides a general framework for such problem but lacks some very crucial details. In this work, we address those crucial points and provide evidence why it should be improved as stated in this study.  In addition, the work of Nguyen and Stepanyan \cite{NguyenStepanyan10} investigates the effect of the engine response time requirements of a generic transport aircraft in severe damage situations associated with the vertical stabilizer. They carried out a study which concentrates on assessing the requirements for engine design for fast responses in an emergency situation. In addition, the use of differential thrust as a propulsion command for the control of directional stability of a damaged transport aircraft was studied by Urnes and Nielsen \cite{UrnesNielsen10} to identify the effect of the change in aircraft performance due to the loss of the vertical stabilizer and to make an improvement in stability utilizing engine thrust as an emergency yaw control mode with feedback from the aircraft motion sensors.  

Existing valuable research in literature provides insight regarding the dynamics of such an extreme scenario, including some unique studies on nonlinear control\cite{Murillo2015} and diverse applications such as adaptive fault detection and isolation methodologies \cite{LiYang2012}. In this paper, in comparison to the existing works, a novel framework and a methodology is provided where Model Reference and Lyapunov based adaptive control methodology is implemented to aid such a damaged aircraft to land safely, with provided asymptotic stability guarantees. 

The paper is organized as the following: The nominal and damaged aircraft models are derived in Section \ref{AC Models}. The plant dynamics of both the nominal (undamaged) aircraft and of the damaged aircraft are investigated in Section \ref{Plant Dynamics}. In Section \ref{Propulsion}, the engine dynamics of the jet aircraft are modeled as a system of differential equations with corresponding time constant and time delay terms to study the engine response characteristic with respect to a differential thrust input. In Section \ref{Differential Thrust}, a novel differential thrust control module is developed to map a rudder input to a differential thrust input. In Section \ref{Open Loop}, the aircraft's open loop system response is investigated. Then in Section \ref{LQR}, a Linear Quadratic Regulator controller is designed to stabilize the damaged aircraft and provide model reference dynamics. In Section \ref{Adaptive}, the Lyapunov stability approach based model reference adaptive control methodology is implemented to test the ability of the damaged aircraft to \emph{mimic} the model aircraft's (reference) response and achieve safe (and stable) operating conditions. In Section \ref{Robustness}, robustness and uncertainty analysis is conducted to test the stability and validate the overall performance of the system in the presence of uncertainty. In Section \ref{Conclusion}, the paper is finalized. 

\section{The Aircraft Models} \label{AC Models}

\subsection{The Nominal Aircraft Model}
For this research, the Boeing 747-100 was chosen as the main application platform to demonstrate that such Lytapunov and Model Reference Adaptive Control (MRAC) based strategy is applicable to commercial aviation jetliners. The data for the nominal (undamaged) Boeing 747-100 are summarized in Table \ref{table:The undamaged aircraft data}.
\begin{table}[htbp!]
\centering
\caption{The nominal (undamaged) aircraft data \cite{747_ac_charac, roskam87}}
\label{table:The undamaged aircraft data}
 \begin{tabular}{|c|c|}
 \hline
 \textbf{Parameter} & \textbf{Value}\\ \hline
%Center of Gravity $(\bar{x}_{cg})$ & 0.25 \\ \hline
$S$ $(ft^2)$ & 5500  \\ \hline
$b$ $(ft)$ & 196  \\ \hline
$\bar{c}$ $(ft)$ & 27.3  \\ \hline
$y_e$ $(ft)$ & 69.83  \\ \hline
$W$ $(lbs)$ &  $6.3663*10^5$\\ \hline
$m$ $(slugs)$ &  $19786.46$\\ \hline 
$I_{xx}$ $(slug*ft^2)$ & $18.2*10^6$  \\ \hline 
$I_{yy}$ $(slug*ft^2)$ & $33.1*10^6$  \\ \hline
$I_{zz}$ $(slug*ft^2)$ & $49.7*10^6$  \\ \hline
$I_{xz}$ $(slug*ft^2)$ & $0.97*10^6$  \\ \hline
${C_L}_{\beta}$	& -0.160 \\ \hline 
${C_L}_{p}$	& -0.340 \\ \hline 
${C_L}_{r}$	& 0.130 \\ \hline 
${C_L}_{\delta_a}$	& 0.013 \\ \hline 
${C_L}_{\delta_r}$	& 0.008 \\ \hline 
${C_N}_{\beta}$	& 0.160 \\ \hline 
${C_N}_{p}$	& -0.026 \\ \hline 
${C_N}_{r}$	& -0.28 \\ \hline 
${C_N}_{\delta_a}$	& 0.0018 \\ \hline 
${C_N}_{\delta_r}$	& -0.100 \\ \hline 
${C_Y}_{\beta}$	& -0.90 \\ \hline 
${C_Y}_{p}$	& 0 \\ \hline 
${C_Y}_{r}$	& 0 \\ \hline 
${C_Y}_{\delta_a}$	& 0 \\ \hline 
${C_Y}_{\delta_r}$	& 0.120 \\ \hline 
\end{tabular}
\end{table} 

%\pagebreak

Taken from Nguyen and Stepanyan \cite{NguyenStepanyan10}, the lateral/directional linear equations of motion of the nominal (undamaged) aircraft, with the intact ailerons and rudder as control inputs, are presented as

\begin{equation}
\left[
\begin{array}{c}
\dot{\phi}\\
\dot{p}\\
\dot{\beta}\\
\dot{r}
\end{array}\right]
= \left[
\begin{array}{cccc}
0 & 1 & 0 & \bar{\theta} \\
0 & L_p & L_{\beta} & L_r \\
\frac{g}{\bar{V}} & \frac{Y_p}{\bar{V}} & \frac{Y_{\beta}+g\bar{\gamma}}{\bar{V}} & \frac{Y_p}{\bar{V}}-1 \\
0 & N_p & N_{\beta} & N_r
\end{array}\right]
\left[
\begin{array}{c}
{\phi}\\
{p}\\
{\beta}\\
{r}
\end{array}\right]
+
\left[
\begin{array}{cc}
0 & 0\\
L_{\delta_a} & L_{\delta_r} \\
\frac{Y_{\delta_a}}{\bar{V}} & \frac{Y_{\delta_r}}{\bar{V}}\\
N_{\delta_a} & N_{\delta_r}
\end{array}\right]
\left[
\begin{array}{c}
{\delta_a}\\
{\delta_r}
\end{array}\right]
\end{equation}
where the states are $\phi$, $p$, $\beta$, and $r$, which represent the roll angle, roll rate, side-slip angle, and yaw rate, respectively. The corresponding control inputs are $\delta_a$ (aileron input/command) and $\delta_r$  (rudder input/command).

\subsection{The Damaged Aircraft Model} \label{damaged AC model}

For the modeling of the damaged aircraft, in case of the loss of the vertical stabilizer, lateral/directional stability derivatives need to be reexamined and recalculated. Since in this study it is assumed that the vertical stabilizer is \emph{completely} lost/damaged, the whole aerodynamic structure will be affected, and the new corresponding stability derivatives have to be calculated, and studied. The lateral/directional dimensionless derivatives that depend on the vertical stabilizer include \cite{nelson98}:

\begin{equation}
 {C_Y}_{\beta}=-\eta \frac{S_v}{S}{{C_L}_\alpha}_v\left(1+\frac{d\sigma}{d\beta}\right)
\end{equation}
 
\begin{equation}
  {C_Y}_r=-2\left(\frac{l_v}{b}\right){{C_Y}_\beta}_{tail}
\end{equation}
   
\begin{equation}   
 {C_N}_{\beta}={{C_N}_{\beta}}_{wf}+\eta_v V_v{{C_L}_\alpha}_v\left(1+\frac{d\sigma}{d\beta}\right)
\end{equation}

\begin{equation}
  {C_N}_r=-2\eta_v V_v\left(\frac{l_v}{b}\right){{C_L}_{\alpha}}_v
\end{equation}

\begin{equation}
  {C_L}_r=\frac{C_L}{4}-2\left(\frac{l_v}{b}\right)\left(\frac{z_v}{b}\right){{C_Y}_\beta}_{tail}
\end{equation}

Due to the loss of the vertical stabilizer, the vertical tail area, volume, and efficiency factor will all be zero; therefore, ${C_Y}_{\beta}={C_Y}_r={C_N}_r=0$. If the vertical stabilizer is assumed to be the primary aerodynamic surface responsible for the weathercock stability, then ${C_N}_{\beta}=0$. Finally, ${C_L}_r=\frac{C_L}{4}$.

In addition, without the vertical stabilizer, the mass and inertia data of the damaged aircraft are going to change, where the values that reflect such a scenario (for the damaged aircraft) are listed in Table \ref{table:The damaged aircraft mass and inertia data}.

\begin{table}[htbp!]
\centering
\caption{The damaged aircraft mass and inertia data}
\label{table:The damaged aircraft mass and inertia data}
 \begin{tabular}{|c|c|}
 \hline
 \textbf{Parameter} & \textbf{Value}\\ \hline
$W$ $(lbs)$ &  $6.2954*10^5$\\ \hline 
$m$ $(slugs)$ &  $19566.10$\\ \hline 
$I_{xx}$ $(slug*ft^2)$ & $17.893*10^6$  \\ \hline 
$I_{yy}$ $(slug*ft^2)$ & $30.925*10^6$  \\ \hline
$I_{zz}$ $(slug*ft^2)$ & $47.352*10^6$  \\ \hline
$I_{xz}$ $(slug*ft^2)$ & $0.3736*10^6$  \\ \hline
\end{tabular}
\end{table} 

In this study, during an event of the loss of the vertical stabilizer, it is proposed that the differential thrust component of aircraft dynamics be utilized as an alternate control input replacing the rudder control to regain stability and control of lateral/directional flight dynamics. Next, the lateral-directional linear equations of motion of the damaged aircraft are presented, with the ailerons $(\delta_a)$, differential thrust $(\delta T)$, and collective thrust $(\Delta T)$ as control inputs \cite{NguyenStepanyan10}, as 

\begin{equation}
\left[
\begin{array}{c}
\dot{\phi}\\
\dot{p}\\
\dot{\beta}\\
\dot{r}
\end{array}\right]
= \left[
\begin{array}{cccc}
0 & 1 & 0 & \bar{\theta} \\
0 & L_p & L_{\beta} & L_r \\
\frac{g}{\bar{V}} & \frac{Y_p}{\bar{V}} & \frac{Y_{\beta}+g\bar{\gamma}}{\bar{V}} & \frac{Y_p}{\bar{V}}-1 \\
0 & N_p & N_{\beta} & N_r
\end{array}\right]
\left[
\begin{array}{c}
{\phi}\\
{p}\\
{\beta}\\
{r}
\end{array}\right]
+
\left[
\begin{array}{ccc}
0 & 0 & 0\\
L_{\delta_a} & \frac{\overline{I_{xz}}y_e}{\overline{I_{xx}}\overline{I_{zz}}-\overline{I_{xz}}^2} & 0\\
\frac{Y_{\delta_a}}{\bar{V}} & 0 & \frac{-\bar{\beta}}{m\bar{V}}\\
N_{\delta_a} & \frac{\overline{I_{xx}}y_e}{\overline{I_{xx}}\overline{I_{zz}}-\overline{I_{xz}}^2} & 0
\end{array}\right]
\left[
\begin{array}{c}
{\delta_a}\\
{\delta T} \\
{\Delta T} \\
\end{array}\right]
\end{equation}

In this case, if the initial trim side-slip angle is zero, then $\Delta T$ does not have any significance in the control effectiveness for a small perturbation around the trim condition \cite{NguyenStepanyan10}, which means that the above equations of motion can be reduced to the final form of governing equations of motion for damaged aircraft as:

\begin{equation}\label{eq:lat_damaged_dynamics}
\left[
\begin{array}{c}
\dot{\phi}\\
\dot{p}\\
\dot{\beta}\\
\dot{r}
\end{array}\right]
= \left[
\begin{array}{cccc}
0 & 1 & 0 & \bar{\theta} \\
0 & L_p & L_{\beta} & L_r \\
\frac{g}{\bar{V}} & \frac{Y_p}{\bar{V}} & \frac{Y_{\beta}+g\bar{\gamma}}{\bar{V}} & \frac{Y_p}{\bar{V}}-1 \\
0 & N_p & N_{\beta} & N_r
\end{array}\right]
\left[
\begin{array}{c}
{\phi}\\
{p}\\
{\beta}\\
{r}
\end{array}\right]
+
\left[
\begin{array}{cc}
0 & 0\\
L_{\delta_a} & \frac{\overline{I_{xz}}y_e}{\overline{I_{xx}}\overline{I_{zz}}-\overline{I_{xz}}^2}\\
\frac{Y_{\delta_a}}{\bar{V}} & 0 \\
N_{\delta_a} & \frac{\overline{I_{xx}}y_e}{\overline{I_{xx}}\overline{I_{zz}}-\overline{I_{xz}}^2}
\end{array}\right]
\left[
\begin{array}{c}
{\delta_a}\\
{\delta T} \\
\end{array}\right]
\end{equation}

\label{eq:lat_damaged_dynamics}

\subsection{Flight Conditions}
In this research, the flight scenario is chosen to be a steady, level cruise flight for the Boeing 747-100 at Mach 0.65 (with the corresponding airspeed of $673[ft/sec]$) at $20,000[ft]$ altitude. It is assumed that at one point during the flight, the vertical stabilizer is completely damaged and the aircraft remains with virtually no vertical stabilizer. 
% The corresponding flight conditions for both the damaged and undamaged aircraft models in this research are summarized in Table \ref{table:Flight conditions} and 

In the followings, the means to control aircraft is investigated in such an extreme case scenario. For this purpose, next nominal (undamaged) and damaged aircraft models are developed for analysis.  

% \begin{table}[htbp!]
% \centering
% \caption{Flight conditions}
% \label{table:Flight conditions}
%  \begin{tabular}{|c|c|}
%  \hline
%  \textbf{Parameter} & \textbf{Value}\\ \hline
%  Altitude $(ft)$ & 20,000 \\ \hline
% Air Density $(slugs/ft^3)$ & 0.001268 \\ \hline
% Airspeed $(ft/s)$ & 673 \\ \hline
% \end{tabular}
% \end{table} 

\section{Plant Dynamics} \label{Plant Dynamics}

With the given flight conditions data and information provided in Table \ref{table:The undamaged aircraft data}, the corresponding state space representation for the lateral/directional dynamics of the nominal (undamaged) Boeing 747-100 are obtained as

\begin{equation}\label{eq:AnBn}
A_n= \left[
\begin{array}{cccc}
0 & 1 & 0 & 0\\
0 & -0.8566 & -2.7681 & 0.3275\\
0.0478 & 0 & -0.1079 & -1 \\
0 & -0.0248 & 1.0460 & -0.2665
\end{array}\right],~~
B_n= \left[
\begin{array}{cc}
0 & 0\\
0.2249 & 0.1384\\
0 & 0.0144\\
0.0118 & -0.6537
\end{array}\right]
\end{equation}

\begin{equation}
 C_n= I_{(4x4)} 
%  \left[
% \begin{array}{cccc}
% 1 & 0 & 0 & 0\\
% 0 & 1 & 0 & 0\\
% 0 & 0 & 1 & 0\\
% 0 & 0 & 0& 1
% \end{array}\right]
,~~
 D_n= 0_{(4x2)} 
%  \left[
% \begin{array}{cc}
% 0 & 0\\
% 0 & 0\\
% 0 & 0\\
% 0 & 0
% \end{array}\right]
\end{equation}

Based on the data for the lateral/directional stability derivatives of the aircraft without its vertical stabilizer (as given in Section II. \ref{damaged AC model}), the lateral/directional representation of the damaged aircraft can be achieved as

\begin{equation} \label{eq:AdBd}
{A_d}= \left[
\begin{array}{cccc}
0 & 1 & 0 & 0\\
0 & -0.8566 & -2.7681 &  0.1008\\
0.0478 & 0 & 0 & -1 \\
0 & -0.0248 & 0 & 0
\end{array}\right],~~
{B_d}= \left[
\begin{array}{cc}
0 & 0\\
0.2249 & 0.0142\\
0 & 0\\
0.0118 & 0.6784
\end{array}\right]
\end{equation}

\begin{equation} \label{eq:Cd}
C_d=C_n= I_{(4x4)} 
% \left[
% \begin{array}{cccc}
% 1 & 0 & 0 & 0\\
% 0 & 1 & 0 & 0\\
% 0 & 0 & 1 & 0\\
% 0 & 0 & 0& 1
% \end{array}\right]
% \end{equation}
% 
% \begin{equation} \label{eq:Dd}
,~~
D_d=D_n= 0_{(4x2)}
% \left[
% \begin{array}{cc}
% 0 & 0\\
% 0 & 0\\
% 0 & 0\\
% 0 & 0
% \end{array}\right]
\end{equation}

Here, ${A_n}$ defines the state matrix of the nominal, undamaged aircraft whereas ${A_d}$ represents the state matrix of the damaged aircraft. Furthermore, ${B_n}$ represents the input matrix where the ailerons $(\delta_a)$ and rudder $(\delta_r)$ are control inputs of the undamaged (nominal) aircraft whereas ${B_d}$ stands for the input matrix of the scenario where the ailerons $(\delta_a)$ and differential thrust $(\delta T)$ are control inputs of the damaged aircraft. It is also worth noting that the structure of the input matrix of the nominal aircraft $(B_n)$ and of damaged aircraft $(B_d)$ remain fairly similar, except for the $b_{32}$ term. In $B_n$, $b_{32}$ equals $0.0144$, but $b_{32}$ equals $0$ in $B_d$, which removes the effect of differential thrust on side-slip angle. 

At this point, we would like to emphasize a very major distinction between the existing studies in literature (including Stepanyan et al. \cite{Stepanyan09}) and this work.

\begin{rmrk}
 It is a very well known (text book) fact that static directional stability (in sense of side-slip) of an aircraft is achieved through the vertical stabilizer. Without the means of the vertical stabilizer, the static directional stability of an aircraft remains very challenging. This is also very clearly shown through the aerodynamic stability derivative analysis (presented in Section-\ref{damaged AC model}) and is evident from the comparison of $B_n$ and $B_d$ as presented in Eq.\ref{eq:AnBn} and Eq.\ref{eq:AdBd}, respectively. If Eq.\ref{eq:AdBd} is observed carefully, it will be clear that in case of the \emph{complete} loss of vertical stabilizer there will be absolutely no correspondence between any of the $\delta_a$ and/or $\delta_T$ inputs vs. the side-slip angle ($\beta$) which is disregarded in the existing works in literature. This creates a unique structure of the input matrix-B and provides challenges in terms of singularities to the control system design. Thus, we propose the following structure of a severely damaged aircraft dynamics to reflect a real-life applicable scenario.
\end{rmrk}

\begin{asmpt}\label{asmpt:bd}
Based on the conducted aerodynamic stability analysis, a damaged aircraft with complete loss of vertical stabilizer preserves the control (input) matrix structure as, 
\begin{equation} \label{eq:B_d}
 B_d = \left[ \begin{array}{cc}0 & 0 \\ b_{21} & b_{22} \\ 0 & 0 \\ b_{41} & b_{42} \end{array}
 \right]
\end{equation}
and is assumed to have limited control authority to represent the damaged vertical stabilizer scenario.
\end{asmpt}

Through Eq.(\ref{eq:B_d}),  Assumption-\ref{asmpt:bd} and based on the geometric analysis of damaged aircraft aerodynamic stability derivatives presented in Section-\ref{damaged AC model}, we emphasize on the fact that, (with the completely damaged vertical fin) there exists mapping \emph{only} to \emph{yaw} and \emph{roll} dynamics, but not necessarily side-slip dynamics. Thus, Assumption-\ref{asmpt:bd} and analysis conducted in Section-\ref{damaged AC model} is important and crucial. 

Next, in addition to the structure of control matrix, here we present the damping characteristics of the nominal and damaged aircraft, and are summarized in Table \ref{table:nominal_damping} and Table \ref{table:damaged_damping}, for further investigation.

\begin{table}[htbp!]
\centering
\caption{Damping characteristics of the nominal aircraft}
\label{table:nominal_damping}
 \begin{tabular}{|c|c|c|c|c|}
 \hline
\textbf{Mode} & \textbf{Pole Location} & \textbf{Damping} & \textbf{Frequency $(1/s)$} & \textbf{Period $(s)$}\\ \hline 
Dutch Roll & $-0.126\pm i1.06$ & $0.118$ & $1.07$ & $5.8822$\\ \hline
Spiral & $-0.0172$ & $1$ & $0.0172$ & $365.2651$\\ \hline
Roll & $-0.963$ & $1$ & $0.963$ & $6.5262$\\
\hline
\end{tabular}
\end{table}

\begin{table}[htbp!]
\centering
\caption{Damping characteristics of the damaged aircraft}
\label{table:damaged_damping}
 \begin{tabular}{|c|c|c|c|c|}
 \hline
\textbf{Mode} & \textbf{Pole Location} & \textbf{Damping} & \textbf{Frequency $(1/s)$} & \textbf{Period $(s)$}\\ \hline 
Dutch Roll & $0.0917\pm i0.43$ & $-0.209$ & $0.439$ & $14.2969$\\ \hline
Spiral & $6.32*10^{-18}$ & $-1$ & $6.32*10^{-18}$ & $9.9486*10^{17}$\\ \hline
Roll & $-1.04$ & $1$ & $1.04$ & $6.0422$\\
\hline
\end{tabular}
\end{table}

\begin{figure}[h!]
  	\centering
 	\includegraphics[width=5truein]{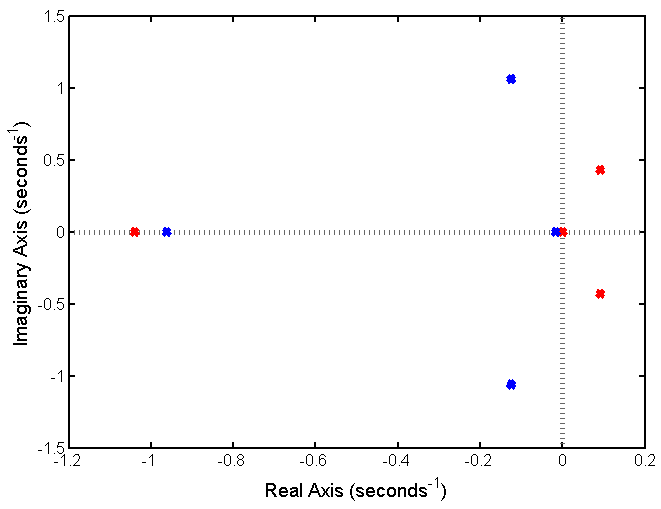}
  	\caption{Pole locations of the nominal \textcolor{blue}{(blue)} and  damaged aircraft \textcolor{red}{(red)}}
 	\label{fig:pole_location}      
 	\end{figure}

Table \ref{table:nominal_damping} shows that all three lateral/directional modes of the nominal aircraft are stable due to the Left Half Plane (LHP) pole locations whereas Table \ref{table:damaged_damping} clearly indicates the unstable nature of the damaged aircraft in the Dutch roll mode by the Right Half Plane (RHP) pole locations. Furthermore, the pole of the spiral mode lies at the origin, which represents very slow (also unstable) dynamics. The only stable mode of the damaged aircraft is the roll mode due to the Left Half Plane (LHP) pole location. The pole locations of both the nominal and damaged aircraft are also illustrated in Fig. \ref{fig:pole_location}.

In the analysis of aircraft dynamics without the vertical stabilizer, propulsion dynamics will play a vital role in the \emph{maneuverability} of the aircraft and deserves a thorough investigation, which is investigated in further detail in the following section.
 	
\section{Propulsion Dynamics} \label{Propulsion}
With emerging advancements in manufacturing processes, structures, and materials, it is a well known fact that aircraft engines have become highly complex systems and include numerous nonlinear processes, which affect the overall performance (and stability) of the aircraft. From the force-balance point of view, this is usually due to the existing coupled and complex dynamics between engine components and their relationships in generating thrust. However, in order to utilize the differential thrust generated by the jet engines as a control input for lateral/directional stability, the dynamics of the engine need to be modeled in order to gain an insight into the response characteristics of the engines.  

Engine response, generally speaking, depends on its time constant and time delay characteristics. Time constant dictates how fast the thrust is generated by the engine, while time delay (which is inversely proportional to the initial thrust level) is due to the lag in engine fluid transport and the inertias of the mechanical systems such as rotors and turbo-machinery blades \cite{NguyenStepanyan10}. 

It is also suggested \cite{NguyenStepanyan10} that the non-linear engine dynamics model can be simplified as a time-delayed second-order linear model as
\begin{equation}
\ddot{T}+2\zeta\omega\dot{T}+\omega^2T= \omega^2T_c(t-t_d)
\end{equation}
where $\zeta$ and $\omega$ are the damping ratio and bandwidth frequency of the closed-loop engine dynamics, respectively; $t_d$ is the time delay factor, and $T_c$ is the thrust command prescribed by the engine throttle resolver angle.

With the time constant defined as the inverse of the bandwidth frequency $(\tau=\frac{1}{\omega})$, and $\zeta$ chosen to be 1 representing a critically damped engine response (to be comparable to existing studies), the engine dynamics can be represented as

\begin{equation}
\left[
\begin{array}{c}
\dot{T}\\
\ddot{T}
\end{array}\right]
= \left[
\begin{array}{cc}
0 & 1\\
\frac{-1}{\tau^2} & \frac{-2}{\tau}
\end{array}\right]
\left[
\begin{array}{c}
T\\
\dot{T}
\end{array}\right]
+
\left[
\begin{array}{c}
0 \\
\frac{1}{\tau^2}
\end{array}\right]
T_c(t-t_d)
\end{equation}

For this study, the Pratt and Whitney JT9D-7A engine is chosen for the application in the Boeing 747-100 example, where the engine itself produces a maximum thrust of 46,500 lbf \cite{boeing747_tech_spec}. At Mach 0.65 and 20,000 feet flight conditions, the engine time constant is 1.25 seconds, and the time delay is 0.4 second \cite{NguyenStepanyan10}.

\begin{figure}[h!]
 \centering
 \includegraphics[width=5truein]{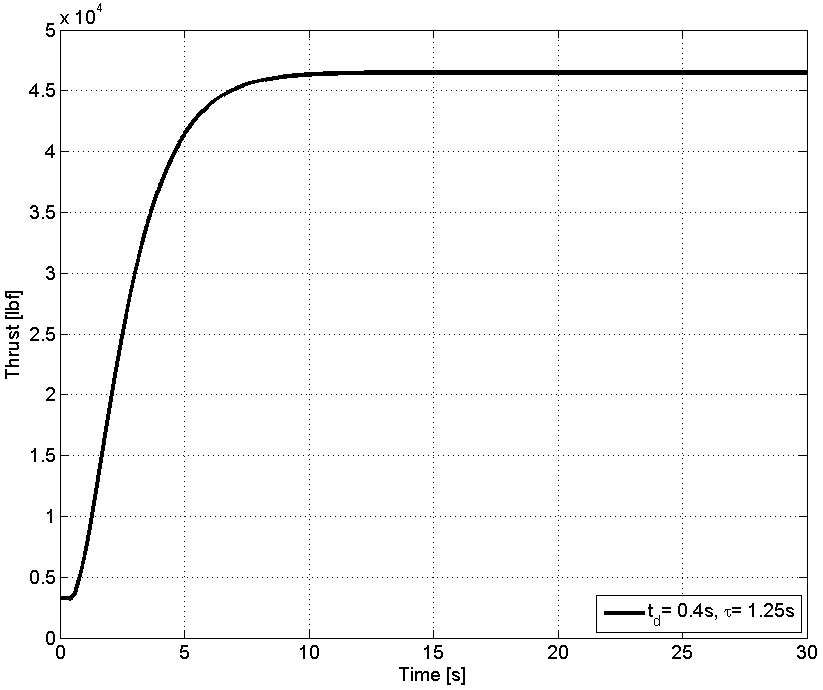}
 \caption{Engine thrust response at Mach 0.65 and 20,000 feet}
 \label{fig:engine_response}
 \end{figure}

The engine thrust response curve at Mach 0.65 and 20,000 feet is, therefore, obtained as shown in Fig. \ref{fig:engine_response}, which provides a useful insight into how the time constant and time delay factors affect the generation of thrust for the JT9D-7A jet engine. At Mach 0.65 and 20,000 feet, with the engine time constant of 1.25 seconds, and the time delay of 0.4 second, it takes approximately ten seconds for the engine to reach steady state and generate its maximum thrust capacity at 46,500 lbf from the trim thrust of 3221 lbf. The increase in thrust generation follows a relatively linear fashion with the engine response characteristic of approximately 12,726 lbf/s during the first two seconds, and then the thrust curve becomes nonlinear until it reaches its steady state at maximum thrust capacity after about ten seconds. This represents one major difference between the rudder and differential thrust as a control input. Due to the lag in engine fluid transport and turbo-machinery inertias, differential thrust (as a control input) cannot respond as instantaneously as the rudder, which has to be taken into account very seriously in control system design.

\section{Differential Thrust as a Control Mechanism} \label{Differential Thrust}

\subsection{Thrust Dynamics and Configuration}
\begin{figure}[htbp!]
 \centering
 \includegraphics[width=6truein]{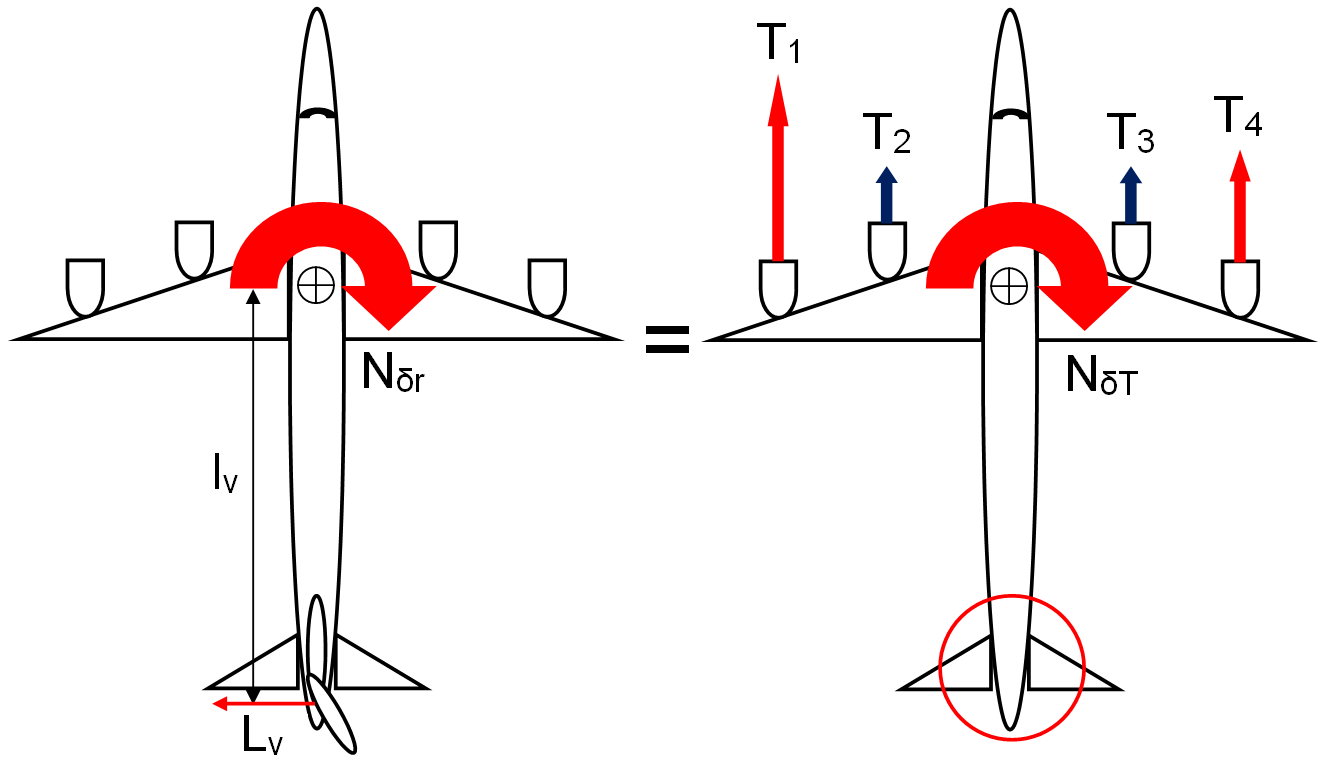}
 \caption{The free body diagram}
 \label{fig:model}
 \end{figure}

In order to utilize differential thrust as a control input for a conventional four-engined (which could be very easily adapted to a twin engine) aircraft, a differential thrust control module is developed to provide a \emph{mapping} between the rudder dynamics and corresponding thrust values. As it is a well known (text book type) concept, the \emph{overall} differential thrust input is defined as the \emph{net} generated thrust. In this scenario, it is defined as the net thrust as a result of engine number 1 and engine number 4 dynamics, while the amounts of thrust generated by remaining engines are kept equal, and fixed as shown in Eqs. (\ref{eq:d_th}-\ref{eq:th23}) to balance the associated torque/moment values. This concept is illustrated in further details in Fig. \ref{fig:model}.

\begin{equation}\label{eq:d_th}
\delta T=T_1-T_4
\end{equation}
\begin{equation}\label{eq:th23}
T_2=T_3
\end{equation}

Engine number 1 and 4 are employed to generate the differential thrust due to the longer moment arm $(y_e)$, which makes the differential thrust more effective as a control for yawing moment. This brings into the picture the need of developing a logic that maps a rudder input to a differential thrust input, which is further explained in the following section.

\subsection{Rudder Input to Differential Thrust Input Mapping Logic}

\begin{figure}[htbp!]
 \centering
 \includegraphics[width=5truein]{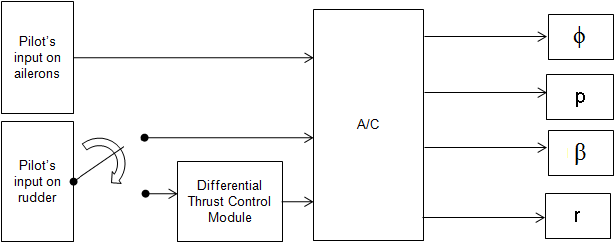}
 \caption{Aircraft control logic diagram}
 \label{fig:control_logics}
 \end{figure}

\begin{figure}[htbp]
 \centering
 \includegraphics[width=5truein]{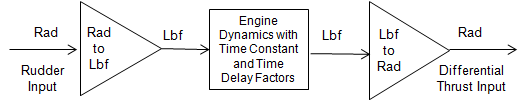}
 \caption{Differential thrust control module}
 \label{fig:dT_control_module}
 \end{figure}

When the vertical stabilizer of an aircraft is intact (i.e. with nominal plant dynamics), ailerons and rudder remain as major control input mechanisms. However, when the vertical stabilizer is damaged, the control effort from the rudder will not respond. To eliminate this mishap, but to still be able to use the rudder demand, a differential thrust control module is introduced in the control logic, as shown in Fig. \ref{fig:control_logics} and Fig. \ref{fig:dT_control_module}, respectively. This differential thrust control module is responsible for mapping corresponding input/output dynamics from the rudder pedals to the aircraft response, so that when the rudder (and whole vertical stabilizer) is completely lost, the rudder input will still be utilized but \emph{switched/converted} to a differential thrust input, which acts as the rudder input for lateral/directional controls. This logic constitutes one of the novel approaches introduced in this paper.

As it can be also seen from Fig. \ref{fig:control_logics} and Fig. \ref{fig:dT_control_module}, the differential thrust control module's function is to convert the rudder (pedal) input to the differential thrust input. In order to achieve that, the rudder (pedal) input (in radians) is mapped to the differential thrust input (in pounds-force) which is then provided into the engine dynamics, as discussed previously in Section \ref{Propulsion}. With this modification, the engine dynamics will dictate how differential thrust is generated, which is then provided as a "virtual rudder" input into the aircraft dynamics. The radian to pound-force conversion, even if it is straight forward text book information, is derived and provided in the next section, for completeness.

\subsection{Radian to Pound-Force Conversion Factor}
Using Fig. \ref{fig:model} and with the steady, level flight assumption, the following relationship can be obtained:
\begin{equation}
N_{\delta_r}= N_{\delta T}
\end{equation}

\begin{equation}
qSb{C_N}_{\delta_r}\delta_r=(\delta T) y_e
\end{equation}
which means the yawing moment by deflecting the rudder and by using differential thrust have to be the same. Therefore, the relationship between the differential thrust control input $(\delta T)$ and the rudder control input $(\delta_r)$ can be obtained as

\begin{equation}
\delta T= \left(\frac{qSb{C_N}_{\delta_r}}{y_e}\right)\delta_r
\end{equation}

Based on the flight conditions at Mach 0.65 and 20,000 feet, and the data for the Boeing 747-100 summarized in Table \ref{table:Flight conditions} and Table \ref{table:The undamaged aircraft data}, the conversion factor for the rudder control input to the differential thrust input is calculated to be 

\begin{equation}
\frac{\delta T}{\delta_r}= -4.43*10^5 \> \frac{lbf}{rad}
\end{equation}

Due to the sign convention of rudder deflection and the free body diagram in Fig. \ref{fig:model}, $\delta_r$ here is negative. Therefore, for the Boeing 747-100, in this study, the conversion factor for the mapping of a rudder input to a differential thrust input is found to be

\begin{equation}
\frac{\delta T}{\delta_r}= 4.43*10^5 \> \frac{lbf}{rad}
\end{equation}

\subsection{Commanded vs. Available Differential Thrust}
% At this point, the worst case scenario is considered, and it is assumed that the aircraft has lost its vertical stabilizer so that the rudder input is converted to the differential thrust input according to the logics discussed previously in this section. 

\begin{figure}[htbp!]
 \centering
 \includegraphics[width=5truein]{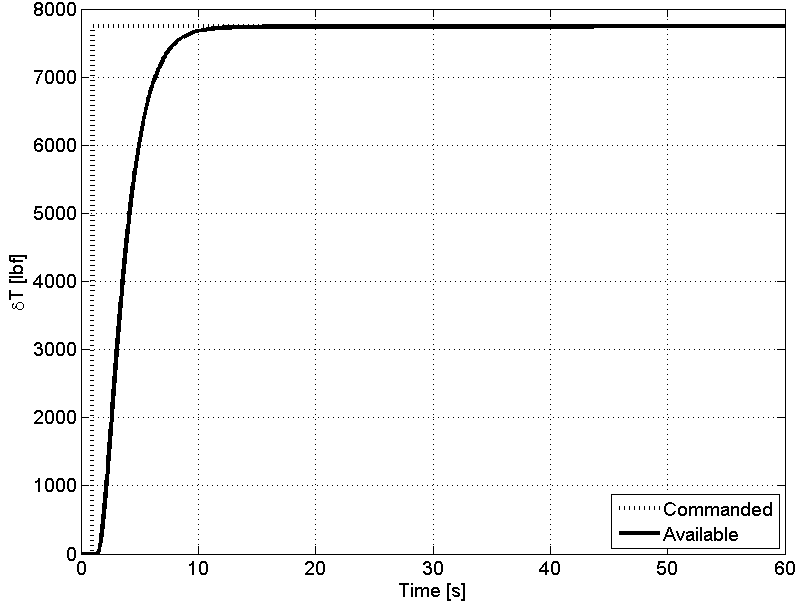}
 \caption{Commanded vs. available differential thrust}
 \label{fig:dT}
 \end{figure}

Unlike the rudder, due to delayed engine dynamics with time constant, there is a major difference in the commanded differential thrust and the available differential thrust as shown in Fig. \ref{fig:dT}.

It can be seen from Fig. \ref{fig:dT} that compared to the commanded differential thrust, the available differential thrust is equal in amount but longer in the time delivery. For a one degree step input on the rudder, the corresponding equivalent commanded and available differential thrust are 7737 lbf,  which are deliverable in ten second duration. Unlike the instantaneous control of the rudder input, there is a lag associated with the use of differential thrust as a control input. This is due to the lag in engine fluid transport and the inertias of the mechanical systems such as rotors and turbo-machinery blades \cite{NguyenStepanyan10}. This is a major design consideration and will be taken into account during the adaptive control system design phase in the following sections. 

\section{Open Loop System Response Analysis}\label{Open Loop}

\begin{figure}[htbp]
 \centering
 \includegraphics[scale=0.8]{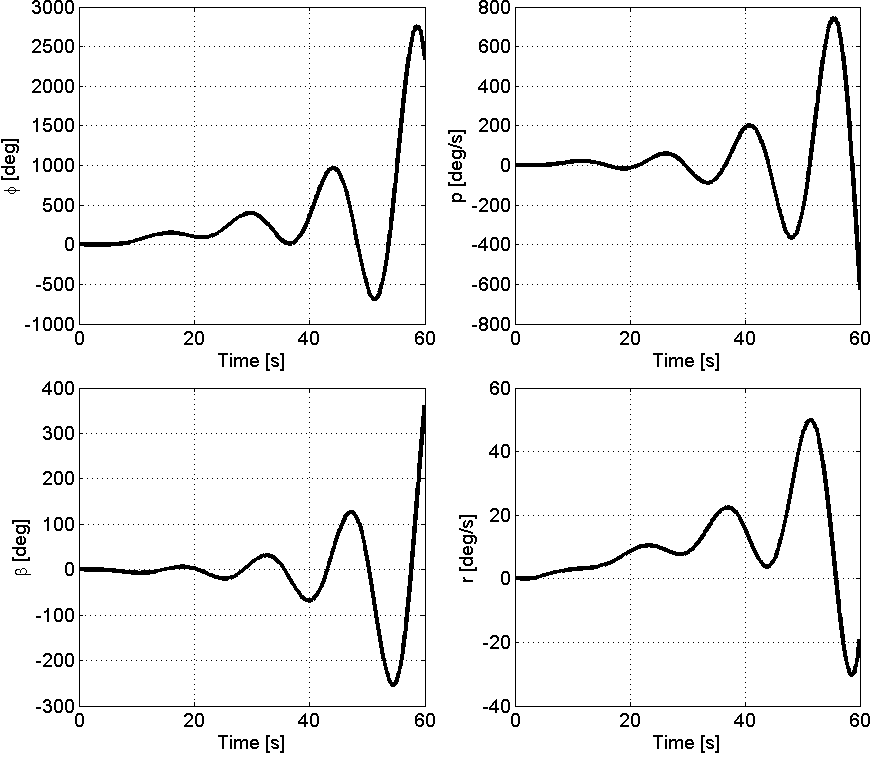}
 \caption{Open loop system response of the damaged aircraft}
 \label{fig:open_loop_response}
 \end{figure}

Following to this, the open loop response characteristics of the aircraft dynamics with a damaged vertical stabilizer to a one degree step input from the ailerons and differential thrust are presented in Fig. \ref{fig:open_loop_response}. It can be clearly seen that when the aircraft is majorly damaged and the vertical stabilizer is lost, the aircraft response to the provided inputs is completely unstable in all four states (as it was also obvious from the pole locations). This means the control authority (or pilot) will not have a chance to stabilize the aircraft in time, which calls for a novel approach to save the damaged aircraft. This is another point where the second novel contribution of the paper is introduced: automatic control strategy to stabilize the aircraft,  which allows safe (i.e. intact) landing of the aircraft.

\section{Linear Quadratic Regulator Design for Model Plant Dynamics}\label{LQR}

\subsection{Background Theory}

As it is well known from literature that optimal control aims for the best (optimal) results within a given set of constraints. An optimal controller is usually designed to minimize a Performance Index (PI), which is generally defined as the ultimate collection of all valuable metrics that are of interest to the designer. On that note, here we present the methodology of obtaining an Linear Quadratic Regulator (LQR) based controller parameters for model plant dynamics that is going to be used in Adaptive Control theory.

Taken from \cite{ogata10}, consider the optimal regulator problem as the following: 

\noindent Given the system equation:
\begin{equation} \label{eq:lqr_ss}
\dot{x}=Ax+Bu
\end{equation}

\noindent Determine the gain matrix $K$ of the optimal control vector:

\begin{equation} \label{eq:lqr_u}
u=-Kx
\end{equation}

\noindent to minimize the Performance Index (PI):

\begin{equation}\label{eq:Jj}
PI=J=\int_{0}^{\infty} (x^TQx+u^TRu)\>dt
\end{equation}

\noindent where Q and R are real, symmetric, positive-definite matrices. It is worth noting that Eq. (\ref{eq:Jj}) represents the Performance Index, in which $x^TQx$ represents the transient energy cost and $u^TRu$ represents the control energy cost. 

\subsection{Stabilizing the Damaged Aircraft with an LQR Controller}\label{sec:lqr}

In this section, we will demonstrate the methodology to obtain closed loop stable model plant dynamics for the damaged aircraft plant. After an iterative process, the state weighting matrix, Q, and the control cost matrix, R, are chosen as

\begin{equation}
Q=10^5 \left[
\begin{array}{cccc}
1 & 0 & 0 & 0\\
0 & 2 & 0 & 0\\
0 & 0 & 0.1 & 0\\
0 & 0 & 0 & 1\\
\end{array}\right]
\end{equation}

\begin{equation}
R=10^3\left[
\begin{array}{cc}
1 & 0\\
0 & 1\\
\end{array}\right]
\end{equation}

where the feedback matrix $K$ is then obtained as

\begin{equation}
K_{LQR}=\left[
\begin{array}{cccc}
9.6697 & 13.2854 & -9.1487 & 0.8729\\
1.9631 & 2.8644 & -12.1067 & 11.5702
\end{array}\right]
\end{equation}

The model plant matrix, $A_m=A_d-B_dK_{LQR}$, then becomes

\begin{equation} \label{eq:Am}
A_m=\left[
\begin{array}{cccc}
 0 & 1 & 0 & 0\\
-2.2026 & -3.8851 & -0.5390 & -0.2595\\
0.0478 & 0 & 0 & -1\\
-1.4455 & -2.1243 & 8.3210 & -7.8597\\
\end{array}\right]
\end{equation}

Next, closed loop response of obtained \emph{model} reference plant is provided. 

\begin{figure}[htbp!]
 \centering
 \includegraphics[scale=0.6]{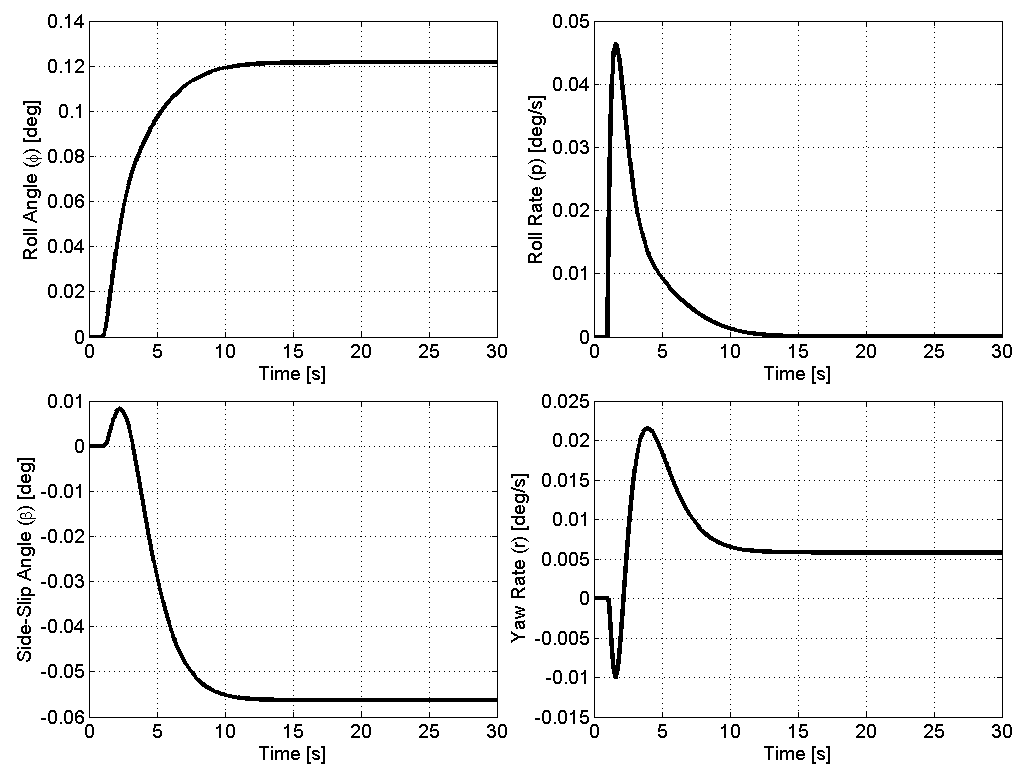}
 \caption{Closed loop response with an LQR controller}
 \label{fig:lqr_state}
 \end{figure}

\begin{figure}[htbp!]
 \centering
 \includegraphics[width=5truein]{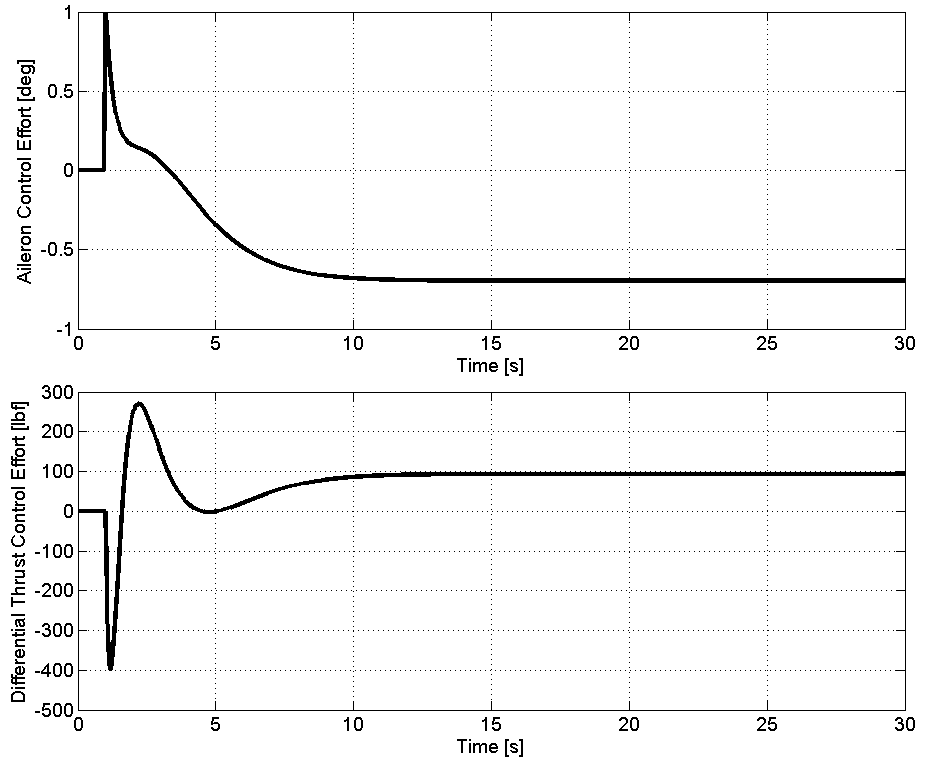}
 \caption{Closed loop control efforts with an LQR controller}
 \label{fig:lqr_control_effort}
 \end{figure}

Compared to the unstable open loop response of the damaged aircraft in Fig. {\ref{fig:open_loop_response} in Section \ref{Open Loop}, the closed loop response is stable in all four states in obtained \emph{model reference plant}: roll angle ($\phi$), roll rate ($p$), side-slip angle ($\beta$), and yaw rate ($r$). In addition, the control efforts for ailerons and differential thrust are also feasible and within the limits of actuator saturation and rate-limiter bounds. From Fig. \ref{fig:lqr_control_effort}, we can also observe that aileron control effort demands the maximum deflection of 1 degree and settles at around -0.7 degree while differential thrust control effort demands a maximum of -400 lbf (negative differential thrust means $T4>T1$) and settles at approximately 100 lbf, which is well within the thrust capability of the JT9D-7A engine.

% In this section, we have proved that it is possible to stabilize the damaged aircraft with a Linear Quadratic Regulator (LQR) controller. Although LQR provides the optimal solution to our problem, it does not guarantee stability in the presence of uncertainty, which is extremely crucial for a damaged aircraft. Therefore, in order to save the damaged aircraft, we will use the LQR controller to stabilize the damaged aircraft (as the inner loop), which will serve as our model aircraft for the adaptive control system design in the next section, and utilize the Lyapunov based model reference adaptive control methodology (as the outer loop) to guarantee stability.

\section{Lyapunov Based Model Reference Adaptive System (MRAS) Controller Design}\label{Adaptive}

To control an aircraft with a fully damaged vertical stabilizer and with no rudder capability can result to be a very stressful and laborious (if not fatal) task for the pilots. This task also requires skills and experience which is hard to possess and execute in extremely stressful moments. In such instances, pilots usually have seconds to react, and as having been witnessed beforehand, coupling between the pilot and unstable aircraft dynamics usually led to a catastrophic crash. Therefore, for the safety of the overall flight, it is crucial for an \emph{online, and adaptive} automatic control system to be developed, tested, and implemented for the aircraft to mitigate accidents and to improve safety, stability, and robustness. As an answer to such need, here, we introduce a novel, Lyapunov stability based adaptive control system design.  

In conventional model reference adaptive control theory, two celebrated and widely used methods are the MIT Rule and the Lyapunov Stability approaches \cite{astrom95}. Because of the Multi-Input-Multi-Output (MIMO) structure of the lateral/directional dynamics, the MIT rule will be left alone due to its relatively weak controllability characteristics in higher order and complex systems \cite{astrom95}. Instead, the powerful nature of Lyapunov based model reference adaptive system (MRAS) contoller design will be utilized.

% In classical adaptive control theory, the Lyapunov stability approach is based on the characteristic of a decreasing kinetic energy function of state dynamics. Because of the reason that kinetic energy of a system is descending, the system is considered approaching its asymptotic stability (equilibrium) point. However, it is a relatively cumbersome task to derive a kinetic energy function for a complex system, but if a candidate function $V(x)$ could be defined, which represents the characteristics of the kinetic energy functions, and if it is descending along the trajectory of the kinetic energy functions, then it can be concluded that the solution of the governing differential equations $\frac{dy}{dx}=f(x)$  will be asymptotically stable. The function $V(x)$ is then called Lyapunov function.

\subsection{Stability Characteristics} \label{adaptive_theorem}

\begin{thm}
For given system dynamics of the damaged aircraft model in Eq. (\ref{eq:AdBd}-\ref{eq:Cd}), there exists a Lyapunov function in form of
\begin{equation} \label{eq:V_x}
V(x)=e^TPe+Tr[(A_d-B_dL-A_m )^T N(A_d-B_dL-A_m )]  
\end{equation}
which guarantees asymptotic stability, if and only if the feedback adjustment law is defined as
\begin{equation}\label{eq:L_dot}
\dot{L}= {({B_d}^TNB_d)}^{-1}{B_d}^TPe{y_d}^T
\end{equation}
\end{thm}

\noindent \emph{Proof:} Let's consider the suggested Lyapunov function taken from \cite{Vempaty11},

\begin{equation} \label{eq:V_x_1}
V(x)=e^TPe+Tr[(A_d-B_dL-A_m)^T N(A_d-B_dL-A_m)]+Tr[(B_dM-B_m)^T R(B_dM-B_m)] 
\end{equation}

For given damaged aircraft dynamics, it is desired that aircraft maintains control (input) matrix structure as defined in \emph{Assumption 1}, leading to $B_m = B_d$, so in Eq. (\ref{eq:V_x_1}),  $Tr[(B_dM-B_m)^T R(B_dM-B_m)]=0$. Therefore,

\begin{equation} \label{eq:V_x_2}
V(x)=e^TPe+Tr[(A_d-B_dL-A_m)^T N(A_d-B_dL-A_m)]  
\end{equation}
Here $N$ is the weighting factor, and $Tr$ is the "Trace" of expression. Also, let's consider $A_m=A_d-B_dK_{LQR}$, as iterated further in Section-\ref{sec:lqr}. It is straight-forward that $V(x)>0,\forall x\neq 0$, $V(0)=0$, and $V(x)$ is continuously differentiable. For given system, error dynamics ($e=y_d-y_m$) becomes
\begin{equation}
\begin{split}
\dot{e}=&~\dot{y_d}-\dot{y_m} \\
=&~(A_dy_d+B_du)-(A_m y_m+B_m u_c)
\end{split}
\end{equation}
With the defined control effort $u=u_c-Ly_d$, 
\begin{equation}
\dot{e}=A_dy_d+B_d(u_c-Ly_d)-A_m y_m-B_m u_c
\end{equation}
where $y_m=y_d-e$. After some algebra, we get
\begin{equation}
\dot{e}=A_m e+(A_d-B_dL-A_m )y_d+(B_d-B_m)u_c
\end{equation}
Again, for given damaged aircraft dynamics, it is desired that aircraft maintains limited control (input) matrix structure as defined in \emph{Assumption 1}, leading to $B_m = B_d$. Thus, 
\begin{equation}
\dot{e}=A_m e+(A_d-B_dL-A_m )y_d
\end{equation}
With $L=L^*+\Delta L $ where  $L^*$ is the constant feedback gain and $\Delta L$ represents the parameter adjustment uncertainty, and $A_d-A_m=B_dL^*$,
\begin{equation}
\dot{e}=A_m e+(A_d-A_m-B_d(L^*+\Delta L))y_d=A_m e+(B_dL^*-B_dL^*-B_d\Delta L)y_d
\end{equation}
Therefore, 
\begin{equation}
\dot{e}=A_m e-B_d\Delta Ly_d
\end{equation}
The derivative of the Lyapunov function from Eq. (\ref{eq:V_x}) can be obtained as
\begin{equation}
\dot{V}(x)=-e^T Qe+2Tr[-{\Delta L}^T B_d^T Pe{y_d}^T+{\Delta L}^T B_d^TNB_d \Delta \dot{L}]
\end{equation}
where ${A_m}^T P+PA_m=-Q$, with positive-definite matrix $Q$ selected to be equal to the Observability Gramian, $C^TC$. From here, it is clear that the negative definite nature of Lyapunov function ($\dot{V}(x)<0$), and therefore, the asymptotic stability of the overall system dynamics is guaranteed when 
\begin{equation}
{\Delta L}^T(-B_d^T Pey^T+ B_d^TNB_d  \Delta \dot{L})=0 
\end{equation}
is satisfied. This leads to the final adaptation law:
\begin{equation}
\dot{L}= {(B_d^TNB_d)}^{-1}B_d^TPe{y_d}^T
\end{equation}
\noindent which guarantees the asymptotic stability $\blacksquare$

\subsection{Simulation Results}

The representative block diagram architecture for the suggested adaptive control system design (based on the Lyapunov stability approach) is illustrated in Fig. \ref{fig:adaptive_diagram}. The ultimate goal of the proposed adaptive control system design is to investigate whether the aircraft with a damaged vertical stabilizer is going to be able to mimic model aircraft dynamics and track the response of the model aircraft or not, by utilizing differential thrust as a control input for lateral/directional dynamics. The control inputs for both plants are one degree step inputs for both the ailerons and differential thrust. It is worth noting that this is an extreme scenario test to see whether the damaged aircraft utilizing differential thrust can hold itself in a continuous yawing and banking maneuver without becoming unstable and losing control.

\begin{figure}[htbp!]
 \centering
 \includegraphics[scale=0.8]{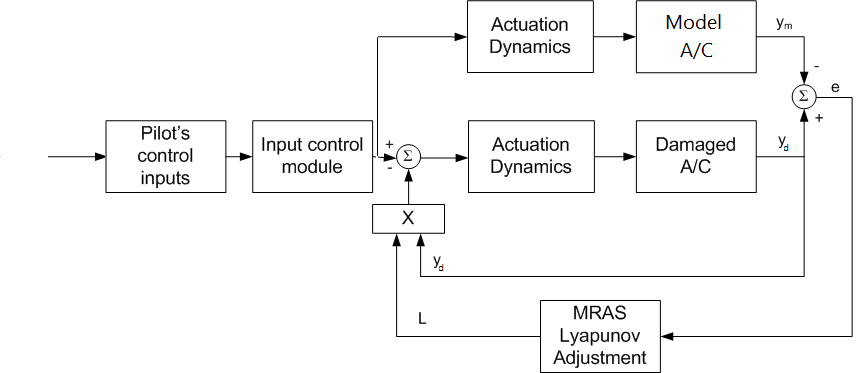}
 \caption{Block diagram structure for adaptive control system}
 \label{fig:adaptive_diagram}
 \end{figure}

As it can also be seen from Fig. \ref{fig:adaptive_diagram}, for both the model and the damaged aircraft, the input signals for the ailerons and rudder are mapped through the input control module, where the rudder input signal is routed through the differential thrust control module and then converted to the differential thrust input following the transformation logics discussed in Section \ref{Differential Thrust} of this paper. 

In order to have a feasible control strategy in real-life situation, limiting factors (such as \emph{saturation limits} and \emph{rate-limiters}) are imposed on the aileron and differential thrust control efforts. The aileron deflection is limited at $\pm$26 degrees \cite{faa_747}. For differential thrust, a differential thrust saturation is set at 43,729 lbf, which is the difference of the maximum thrust and trimmed thrust values of the JT9D-7A engine. In addition, a rate limiter is also imposed on the thrust response characteristic at 12,726 lbf/s as discussed in Section \ref{Propulsion}.

Following to that, the simulation results of the adaptive control system model are presented in Fig. \ref{fig:adaptive_state}. As shown in Fig. \ref{fig:adaptive_state}, after only 15 seconds, all four states of the aircraft's lateral/directional dynamics reach steady state values. It can also be clearly seen that after a time interval of 15 seconds the damaged aircraft plant can mimic the model aircraft plant where the errors are minimized as shown in Fig. \ref{fig:adaptive_error}. This demonstrates the functionality of the Lyapunov based adaptive control system design in such an extreme scenario.
 
\begin{figure}[htbp!]
 \centering
 \includegraphics[scale=0.6]{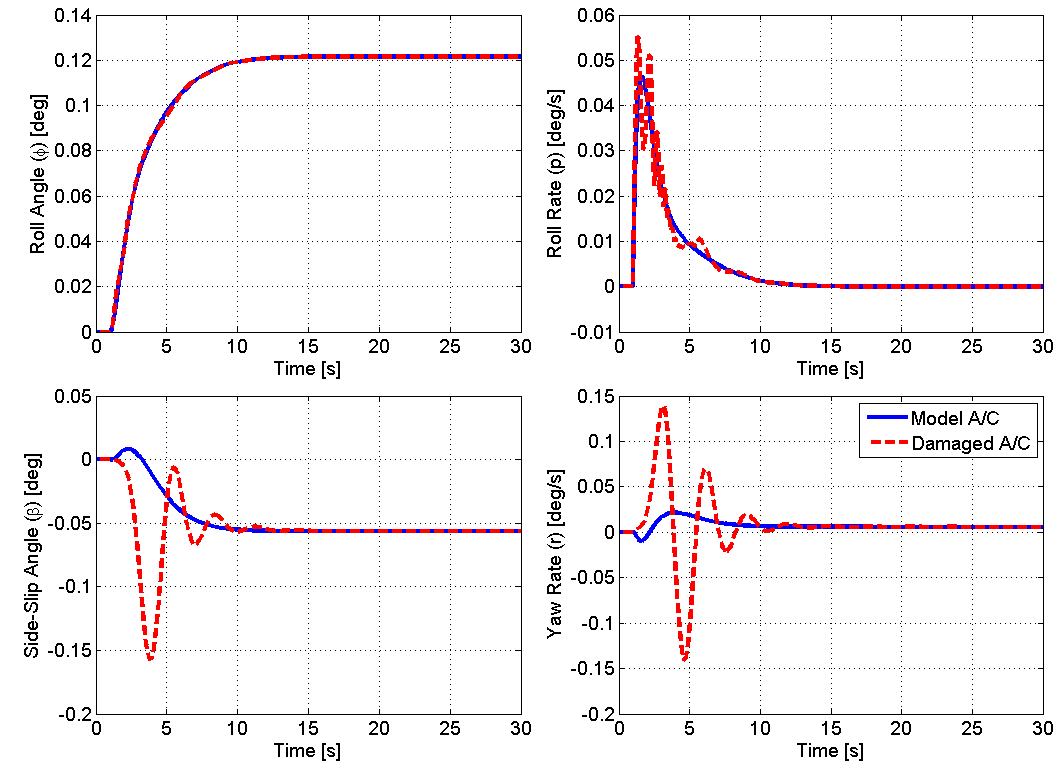}
 \caption{Adaptive control outputs}
 \label{fig:adaptive_state}
 \end{figure}
 
\begin{figure}[htbp!]
 \centering
 \includegraphics[scale=0.6]{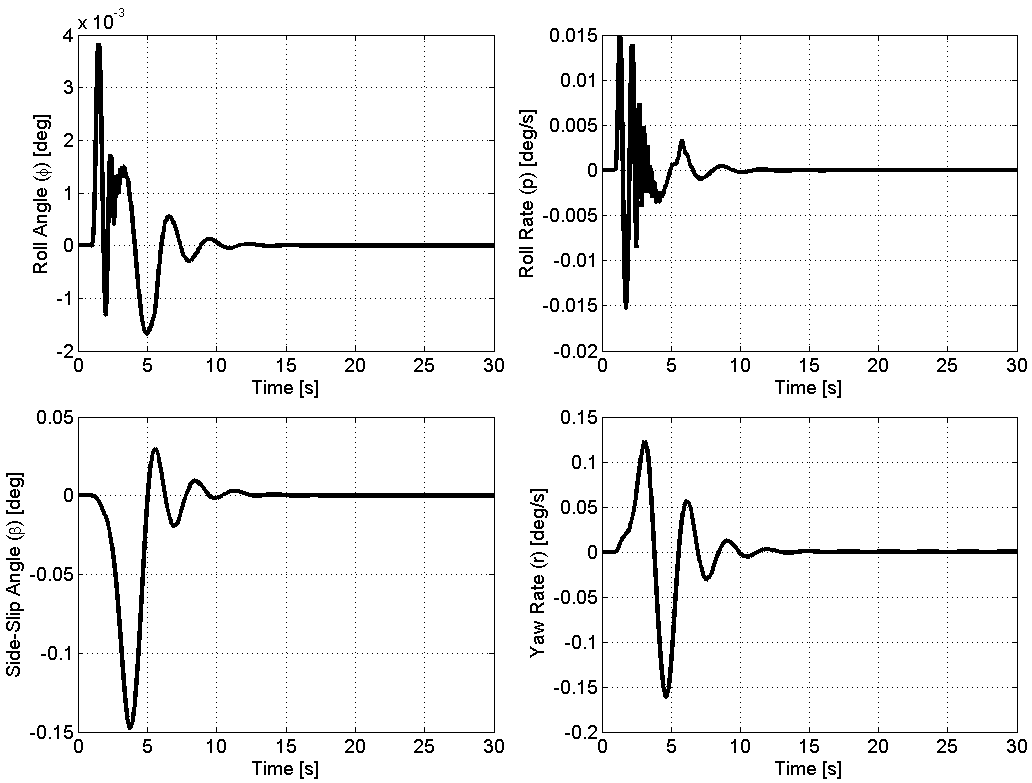}
 \caption{Adaptive error signals}
 \label{fig:adaptive_error}
 \end{figure}

From Fig. \ref{fig:adaptive_error}, it can be observed that the error signals for all four lateral/directional states are diminished after 15 seconds. However, this comes at the cost of slightly higher control effort demand, as shown in Fig. \ref{fig:adaptive_control_effort}, which are still within control limits and without any saturation of the actuators.

\begin{figure}[htbp!]
 \centering
 \includegraphics[width=5truein]{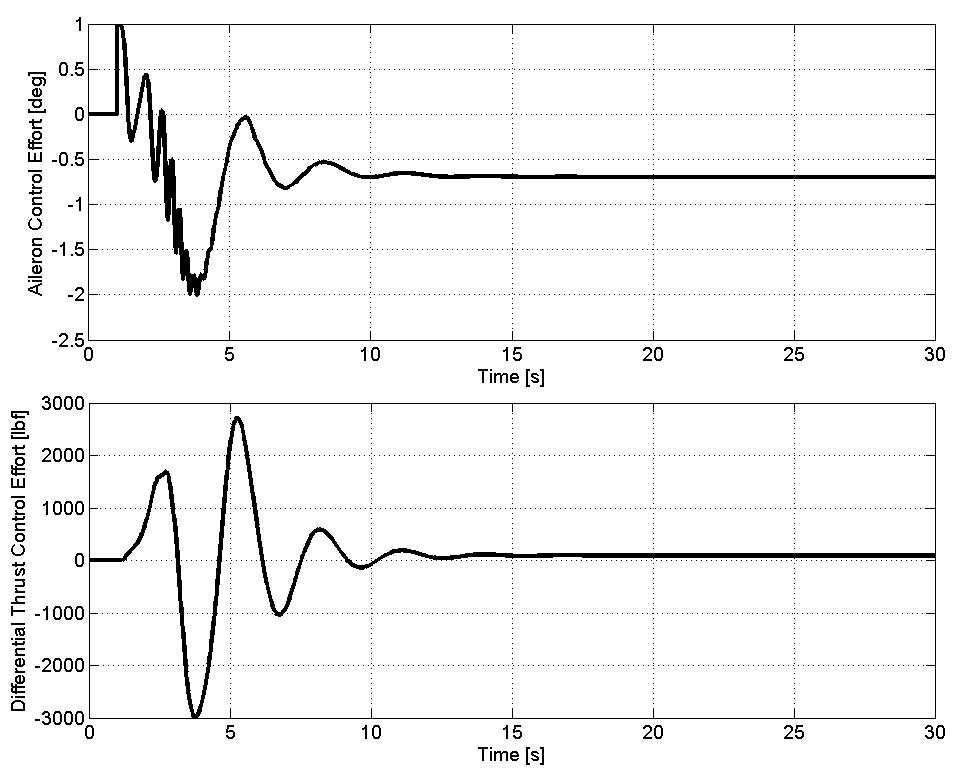}
 \caption{Adaptive control efforts}
 \label{fig:adaptive_control_effort}
 \end{figure}

The aileron control effort, as indicated by Fig.\ref{fig:adaptive_control_effort}, calls for the maximum deflection of about -2 degrees and reaches steady state at approximately -0.7 degree of deflection after 15 seconds responding to a one degree step input. This aileron control effort is very reasonable and achievable if the ailerons are assumed to have instantaneous response characteristics by neglecting the lag from actuators or hydraulic systems. The differential thrust control effort demands a maximum differential thrust of -3000 lbf (negative differential thrust means $T4>T1$), which is within the thrust capability of the JT9D-7A engine, and the differential thrust control effort reaches steady state at around 85 lbf  after 15 seconds. Therefore, it can be concluded that the adaptive control system design with the utilization of differential thrust as a control input is proven to save the damaged aircraft by making it behave like the model aircraft, but the feasibility of the adaptive control method depends heavily on the thrust response characteristics of the aircraft jet engines.

\section{Robustness Analysis of Adaptive Control System Design}\label{Robustness}
The robustness of the adaptive system design presented in this paper is investigated by the introduction of 30\% of full block, additive uncertainty into the plant dynamics of the damaged aircraft, to test its ability to track the reference response of the model aircraft in the presence of uncertainty. Fig. \ref{fig:adaptive_diagram_U} shows the logic behind the adaptive control system design in the presence of uncertainty.

\begin{figure}[htbp!]
 \centering
 \includegraphics[scale=0.8]{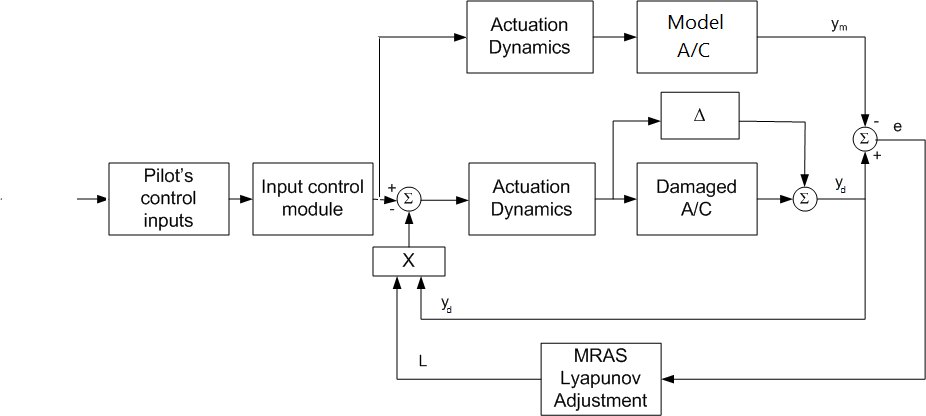}
 \caption{Block diagram for adaptive control system in the presence of uncertainty}
 \label{fig:adaptive_diagram_U}
 \end{figure}
 
One thousand Monte-Carlo simulations were conducted to test the robustness of the damaged plant in the presence of uncertainty. The state responses in the presence of 30\% uncertainty are shown in Fig. \ref{fig:adaptive_state_U}. It is obvious that the adaptive control system design is able to perform well under given uncertain conditions and the damaged aircraft can follow/mimic the response of the model aircraft only after approximately 15 seconds. In that sense, the uncertain plant dynamics are well within the expected bounds.
 
 \begin{figure}[htbp!]
 \centering
 \includegraphics[scale=0.6]{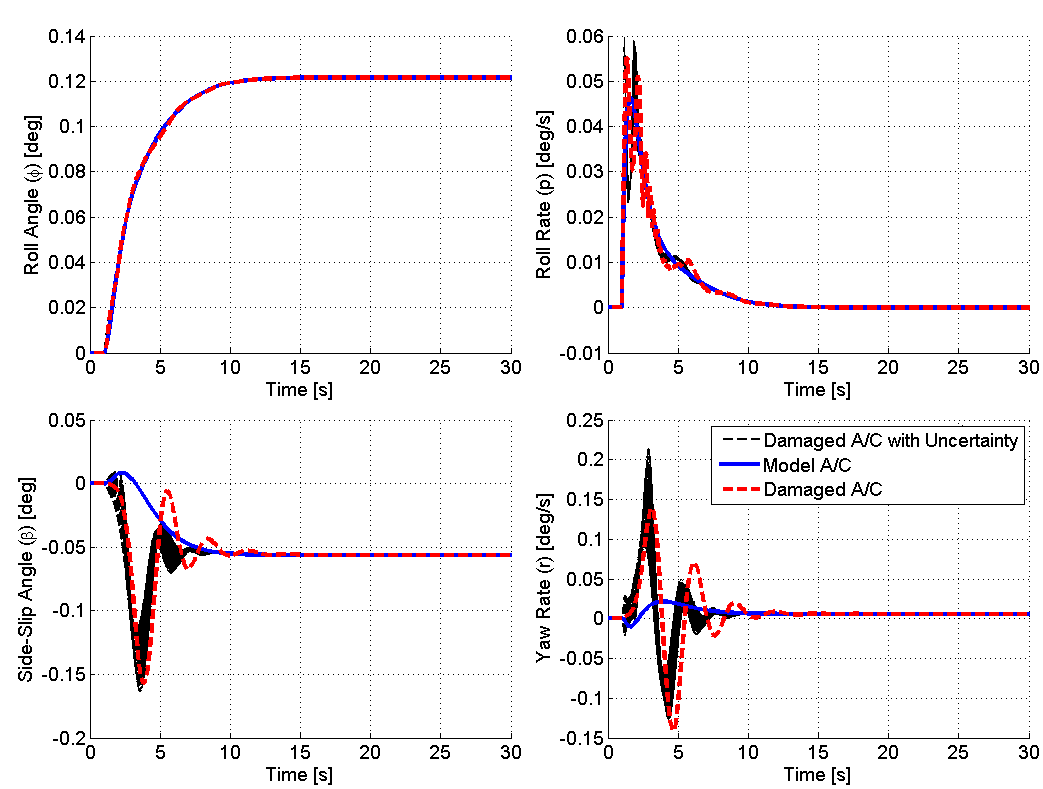}
 \caption{Adaptive control outputs in the presence of 30\% uncertainty}
 \label{fig:adaptive_state_U}
 \end{figure}
 
The robustness of the adaptive control system design can be further illustrated in Fig. \ref{fig:adaptive_error_U} that all the error signals reach steady state and converge to zero only after 15 seconds. However, these favorable characteristics come at the expense of the control effort from the ailerons and differential thrust as shown in Fig. \ref{fig:adaptive_control_effort_U}.

\begin{figure}[htbp!]
 \centering
 \includegraphics[scale=0.6]{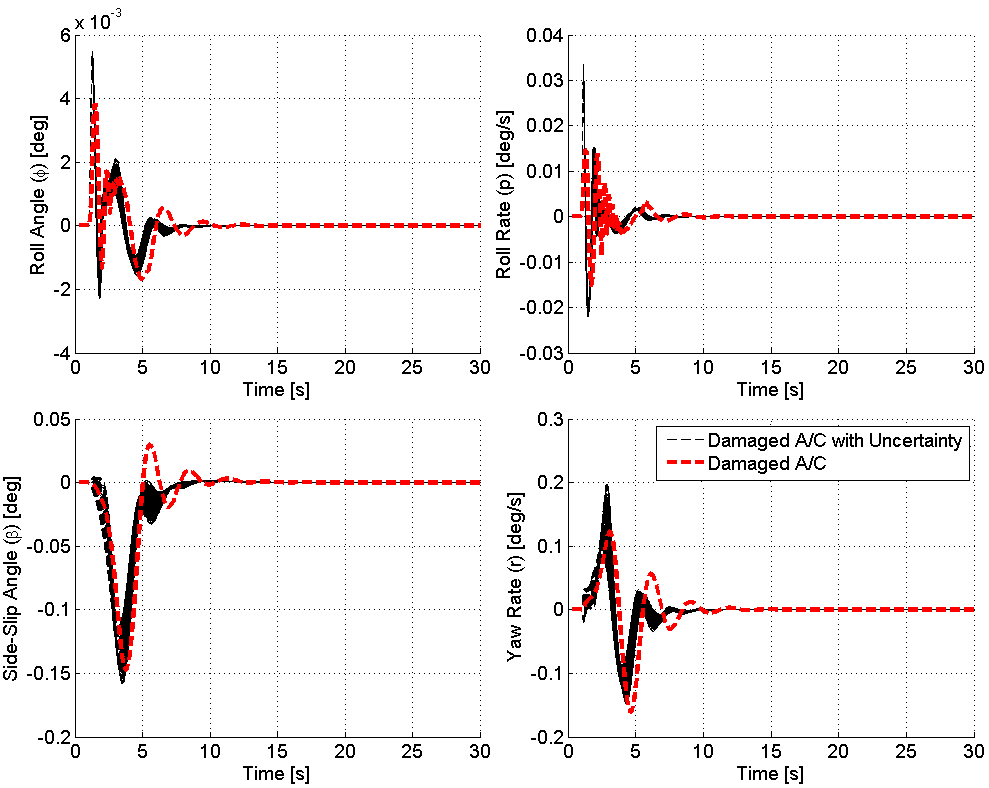}
 \caption{Adaptive error signals in the presence of 30\% uncertainty}
 \label{fig:adaptive_error_U}
 \end{figure}
  
 \begin{figure}[htbp!]
 \centering
 \includegraphics[width=5truein]{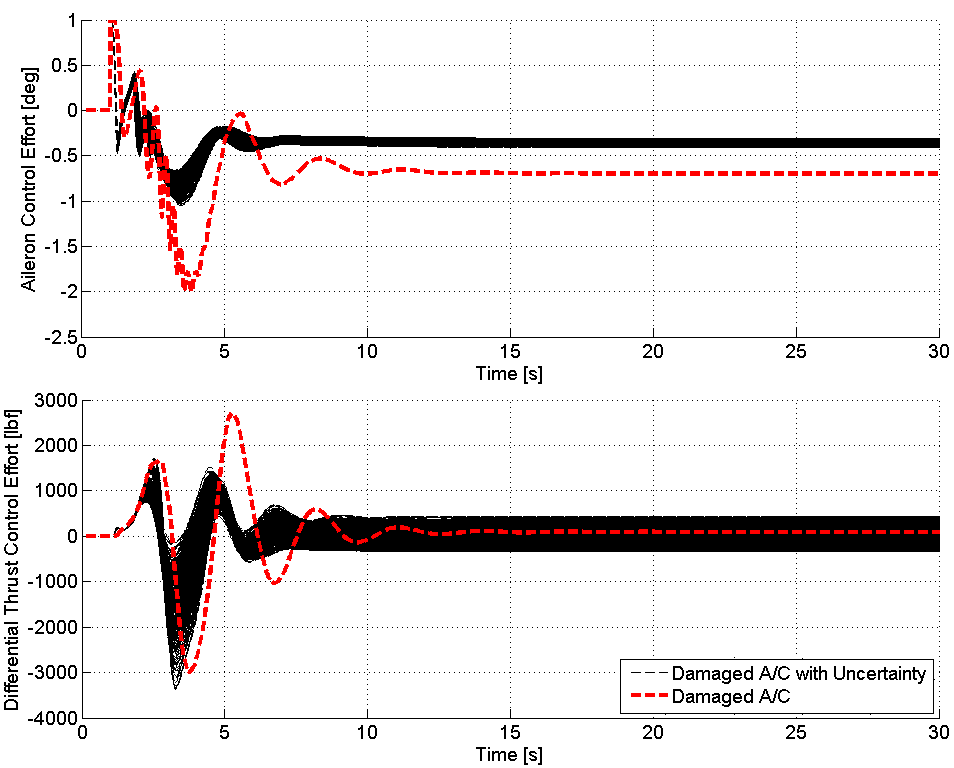}
 \caption{Adaptive control efforts in the presence of 30\% uncertainty}
 \label{fig:adaptive_control_effort_U}
 \end{figure}

According to Fig. \ref{fig:adaptive_control_effort_U}, when there is 30\% full block, additive uncertainty, the aileron control demands the maximum deflection of approximately -1 degree and reaches steady state at around -0.4 to -0.3 degree after 15 seconds. The aileron control effort demands are reasonable and feasible due to the limiting factor of $\pm$ 26 degrees of the aileron deflection \cite{faa_747} and the assumption that ailerons have instantaneous response characteristics by neglecting the lag from actuators or hydraulic systems.
 
As for differential thrust, when there is 30\% uncertainty, the differential thrust control demands at maximum approximately -3400 lbf (negative differential thrust means $T4>T1$), which is within the thrust capability of the JT9D-7A engine, and the differential thrust control effort reaches steady state at around the range of -350 lbf to 450 lbf  after 15 seconds. Again, due to the differential thrust saturation set at 43,729 lbf and the thrust response limiter set at 12,726 lbf/s, this control effort of differential thrust in the presence of uncertainty is achievable in real life situation. 
  
\section{Conclusion} \label{Conclusion}
This paper studied the utilization of differential thrust as a control input to help a Boeing 747-100 aircraft with a damaged vertical stabilizer to regain its lateral/directional stability. 

Throughout this paper, the necessary nominal and damaged aircraft models were constructed, where lateral/directional equations of motion were revisited to incorporate differential thrust as a control input for the damaged aircraft. Then the plant dynamics of both the nominal (undamaged) aircraft and of the damaged aircraft were investigated, and a special aerodynamic-stability-derivative-case (due to severely damaged aircraft geometry) is derived. The engine dynamics of the jet aircraft was modeled as a system of differential equations with engine time constant and time delay terms to study the engine response time with respect to a commanded thrust input. Next, the novel differential thrust control module was then presented to map a rudder input to a differential thrust input. The Linear Quadratic Regulator controller was designed to provide a stabilized model for the damaged aircraft dynamics. The ability of the damaged aircraft to track and mimic the behavior of the model aircraft in an extreme scenario was illustrated through the adaptive control system design based on the Lyapunov stability theory. Demonstrated results showed that the unstable open-loop damaged plant dynamics could be stabilized using adaptive control based differential thrust methodology. Further analysis on robustness showed that uncertain plant dynamics was able to follow model plant dynamics with asymptotic stability, in the presence of 30\% full block, additive uncertainty, associated with damaged aircraft dynamics. 

Overall, this framework provides an automatic control methodology to save the severely damaged aircraft and avoid the dangerous coupling of the aircraft and pilots, which led to crashes in a great number of commercial airline incidents. Furthermore, it has also been concluded that due to the heavy dependence of the differential thrust generation on the engine response, in order to better incorporate the differential thrust as an effective control input in a life-saving scenario, major developments in engine response characteristics are also desired to better assist such algorithm.

\section*{Acknowledgments}
The authors of this paper would like to thank Dr. Ping Hsu from the Department of Electrical Engineering and Dr. Fei Sun from the Flight Control Systems (FCS) and UAV Laboratory at San Jos\'{e} State University for the fruitful discussions and valuable feedback.

%\section{References}

\end{document}